\documentclass[10pt,twoside]{article}
\usepackage{Latex-document, amssymb, latexsym, amsmath}

\renewcommand{\thesection}{\arabic{section}}
\renewcommand{\thesubsection}{\thesection.\arabic{subsection}.}

\def\addsec{\addtocounter{section}{1} \setcounter{theorem}{0}}

\title{\bf  Hard Constraints and the Bethe Lattice: \vskip -2mm
Adventures at the Interface of \vskip -2mm
Combinatorics and Statistical Physics\vskip 6mm}

\author{Graham R. Brightwell\thanks{Department of Mathematics, London School of Economics,
Houghton St., London WC2A 2AE England. E-mail:
g.r.brightwell@lse.ac.uk} \quad {\bf Peter
Winkler}\vspace*{-0.5cm}\thanks{Bell Labs 2C-365, 700 Mountain
Ave., Murray Hill NJ 07974-0636, USA. E-mail: pw@lucent.com}}
\date{\vspace{-8mm}}

\newtheorem{theorem}{Theorem}[section]

\newtheorem{definition}[theorem]{Definition}

\newtheorem{conjecture}[theorem]{Conjecture}

\def\Z{{\mathbb Z}}

\def\T{{\mathbb T}}

\newcommand{\restr}{\! \upharpoonright \!}

\newcommand{\vhi}{{\varphi}}

\renewcommand{\hom}{{\rm Hom}}

\newcommand{\given}{\; \rule[-1mm]{.1mm}{4mm} \;}

\begin{document}

\maketitle

\thispagestyle{first} \setcounter{page}{605}

\begin{abstract}\vskip 3mm

Statistical physics models with hard constraints, such as the discrete
hard-core gas model (random independent sets in a graph), are inherently
combinatorial and present the discrete mathematician with a relatively
comfortable setting for the study of phase transition.

In this paper we survey recent work (concentrating on joint work of the
authors) in which hard-constraint systems are modeled by the space
$\hom(G,H)$ of homomorphisms from an infinite graph $G$ to a fixed
finite constraint graph $H$.  These spaces become sufficiently tractable
when $G$ is a regular tree (often called a Cayley tree or Bethe lattice)
to permit characterization of the constraint graphs $H$ which admit
multiple invariant Gibbs measures.

Applications to a physics problem (multiple critical points for
symmetry-breaking) and a combinatorics problem (random coloring), as
well as some new combinatorial notions, will be presented.
\vskip 4.5mm


\end{abstract}

\vskip 12mm

\section*{1. Introduction} \addsec

\vskip-5mm \hspace{5mm}

Recent years have seen an explosion of activity at the interface
of graph theory and statistical physics, with probabilistic
combinatorics and the theory of computing as major catalysts.
The concept of ``phase transition'', which a short time ago
most graph theorists would barely recognize, has now appeared
and reappeared in journals as far from physics as the
{\em Journal of Combinatorial Theory} (Series B).

Traffic between graph theory and statistical physics is already
heavy enough to make a complete survey a book-length proposition,
even if one were to assume a readership with knowledge of both fields.

This article is intended for a general mathematical audience, {\em not}
necessarily acquainted with statistical physics, but it is not
to serve as an introduction to the field.  Readers are referred
to texts such as \cite{B, G, PF,R} for more background.  We will present
only a small part (but we hope an interesting one) of the interface
between combinatorics and statistical physics, with just enough
background in each to make sense of the text.  We will focus on the most
combinatorial of physical models---those with hard constraints---and
inevitably on the authors' own research and related work.

We hope it will be clear from our development that there is
an enormous amount of fascinating mathematics to be uncovered
by studying statistical physics, quite a lot of which has been or
will be connected to graph theory.  What follows is only a sample.

\section*{2. Random independent sets} \addsec

\vskip-5mm \hspace{5mm}

In what follows a {\em graph} $G = \langle V,E \rangle$ consists of
a set $V$ (finite or countably infinite) of {\em nodes} together
with a set $E$ of {\em edges}, each of which is an unordered
pair of nodes.  We will sometimes permit loops (edges of the
form $\{v,v\}$) but multiple edges will not be considered or needed.
We write $u \sim v$, and say that $u$ is ``adjacent'' to $v$,
if $\{u,v\} \in E$; a set $U \subset V$ is said to be
{\em independent} if it contains no edges.

The {\em degree} of a node $u$ of $G$ is the number of nodes
adjacent to $u$; all graphs considered here will be {\em locally finite},
meaning that all nodes have finite degree.  A {\em path} in $G$
(of length $k$) is a sequence $u_0,u_1, \dots, u_k$ of distinct
nodes with $u_i \sim u_{i+1}$; if in addition $u_k \sim u_0$ we
have a {\em cycle} of length $k\!+\!1$.  If every two nodes of $G$ are
connected by a path, $G$ is said to be {\em connected}.

The plane grid $\Z^2$ is given a graph structure by putting
$(i,j) \sim (i',j')$ iff $|i\!-\!i'|+|j\!-\!j'|=1$.
Let us carve out a big piece of $\Z^2$, say
the box $B_n^2 := \{(i,j) \given -n \le i,j \le n\}$.  Let
$I$ be a uniformly random independent set in the graph $B_n^2$;
in other words, of all sets of nodes (including the empty set)
not containing an edge, choose one uniformly at random.
What does it look like?

Plate 1 shows such an $I$ (in a rectangular region).  Here the
nodes are represented by squares, two being adjacent if they have
a common vertical or horizontal border segment.  The sites belonging
to $I$ are colored, the others omitted.  It is
by no means obvious how to obtain such a random independent
set in practice; one cannot simply choose points one at a
time subject to the independence constraint.  In fact the set
in Plate 1 was generated by {\em Markov chain mixing}, an
important and fascinating method in the theory of computing,
which has by itself motivated much recent work at the
physics-combinatorics interface.

It is common in statistical physics to call the nodes of $B_n^2$
{\em sites} (and the edges, {\em bonds}).  Sites in $I$ are
said to be {\em occupied}; one may imagine that each occupied site
contains a molecule of some gas, any two of which must be at
distance greater than 1.

In the figure, {\em even} occupied sites (nodes $(i,j) \in I$
for which $i\!+\!j \equiv 0 \mod 2$) are indicated by one color,
odd sites by another.  A certain tendency for colors to clump
may be observed; understandably, since occupied sites of the
same parity may be as close as $\sqrt 2$ (in the Euclidean norm)
but opposite-parity particles must be at least $\sqrt 5$ apart.

It stands to reason that if more ``particles'' were forced into
$I$, then we might see more clumping.  Let us weight the
independent sets according to size, as follows: a positive
real $\lambda$, called the {\em activity} (or sometimes
{\em fugacity}) is fixed, and then each independent set $I$
is chosen with probability proportional to $\lambda^{|I|}$.
We call this the ``$\lambda$-measure''.
Of course if $\lambda=1$ we are back to the uniform measure,
but if $\lambda > 1$ then larger independent sets are favored.

The $\lambda$-measure for $\lambda \not= 1$ is, to a physicist,
no less natural than the uniform.  In a combinatorial setting,
such a measure might arise e.g.\ if the particles happen to
be of two different types, with all ``typed'' independent
sets equiprobable; then the probability that a particular set
of sites is occupied is given by the $\lambda$-measure with $\lambda=2$.

Plate 2 shows a random $I$ chosen when $\lambda=3.787$.  The
clusters have grown hugely as more particles were packed in.
Push $\lambda$ up just a bit more, to $3.792$, and something like
Plate 3 is the result: one color (parity) has taken over,
leaving only occasional islands of the other.

Something qualitative has changed here, but what exactly?
The random independent sets we have been looking at constitute
what the physicists call the {\em hard-core lattice gas model},
or ``hard-core model'' for short.  Readers are referred to the
exceptionally readable article \cite{BS1} in which many nice
results are obtained for this model\footnote{Readers, however,
are cautioned regarding conducting a web search with key-words
``hard-core'' and ``model''.}.  On the plane grid, the
hard-core model has a ``critical point'' at activity about 3.79, above
which the model is said to have experienced a {\em phase transition}.

\section*{3. What is phase transition?} \addsec

\vskip-5mm \hspace{5mm}

There is no uniformity even among statistical physicists regarding
the definition of phase transition; in fact, there is even
disagreement about whether the ``phases'' above are the even-dominated
versus odd-dominated configurations at high $\lambda$, or the
high-$\lambda$ regime versus the low.  Technical definitions
involving points of non-analyticity of some function miss the
point for us.

The point really is that a slight change in a parameter governing
the local behavior of some statistical system, like the hard-core model,
can produce a global change in the system, which may be evidenced
in many ways.  For example, suppose we sampled many independent
sets in $B_n^2$ at some fixed $\lambda$, and for each computed
the ratio of the number of even occupied sites to the number of
odd.  For low $\lambda$ these numbers would cluster around $1/2$,
but for high $\lambda$ they would follow a bimodal distribution; and the
larger the box size $n$, the sharper the transition.

Here's another, more general, consideration.  Suppose we look only at
independent sets which contain all the even sites on the boundary
of $B_n^2$.  For these $I$ the origin would be more likely to be
occupied than, say, one of its odd neighbors.   As $n$ grows, this
``boundary influence'' will fade---provided $\lambda$ is low.  But
when $\lambda$ is above the critical point, the boundary values tend
to make $I$ an even-dominated set, giving any even site, no
matter how far away from the boundary, a non-disappearing advantage
over any odd one.

Computationally-minded readers might be interested in a third
approach.  Suppose we start with a fixed independent set $I_0$,
namely the set of all even sites in $B_n^2$, and change it one
site at a time as follows: at each tick of a clock we choose a site
$u$ at random.  If any of $u$'s neighbors is occupied, we do nothing.
Otherwise we flip a biased coin and with probability
$\lambda/(1+\lambda)$ we put $u$ in $I$ (where it may already have
been), and with probability $1/(1+\lambda)$ we remove it (or leave
it out).  The result is a Markov chain whose states are independent
sets and whose stationary distribution, one can easily verify, is
exactly our $\lambda$-measure.  Thus if we do this for many
steps, we will have a nearly perfect sample from this distribution---but
how many steps will that take? We {\em believe} that when $\lambda$
is below its critical value, only polynomially (in $n$) steps are
required---the Markov chain is said to be {\em rapidly mixing};
even polylogarithmic, if we count the number of steps per site.  But for
high $\lambda$ it appears to take time exponential in $n$ (or
perhaps in $\sqrt{n}$) before we can expect to see an odd-dominated
independent set.  The exact relationship between phase transition
and Markov chain mixing is complex and the subject of much study.

All these measures rely on taking limits as the finite box $B_n^2$
grows; the very nice discovery of Dobrushin, Lanford and Ruelle
\cite {D, LR} is that
there is a way to understand the phenomenon of phase transition
as a property of the infinite plane grid.  The idea is to extend
the $\lambda$-measure to a probability distribution on independent sets
on the whole grid, then ask whether the extension is unique.

We cannot extend the definition of the $\lambda$-measure
directly since $\lambda^{|I|}$ is generally infinite, but we
can ask that it behave locally like the finite measure.
We say that a probability distribution $\mu$ on independent sets
in the plane grid is a {\em Gibbs measure} if for any site
$u$ the probability that $u$ is in $I$, given the sites in
$I \cap (\Z^2 \setminus \{u\})$, is $\lambda/(1+\lambda)$ if the
neighborhood of $u$ is unoccupied and, of course, 0 otherwise.

It turns out that Gibbs measures always exist (here, and in
far greater generality) but may or may not be unique.  When
there is more than one Gibbs measure we will say that there is
a phase transition.  For the hard-core model on $\Z^2$, there
is a unique Gibbs measure for low $\lambda$; but above the
critical value, there is a Gibbs measure in which the even
occupied sites are dominant and another in which the odd sites are
dominant (all other Gibbs measures are convex combinations of these
two).  How can you construct these measures?  Well, for example,
the even measure can be obtained as a limit of $\lambda$-measures
on boxes whose even boundary sites are forced to be in $I$.
The fact that the boundary influence does not fade (in the
high $\lambda$ case) implies that the even and odd Gibbs
measures are different.

We have noted that the critical value of $\lambda$ for the
hard-core model on $\Z^2$ is around 3.79.  This is an
empirical result and all we mathematicians can prove is
that there is at least one critical point, and all such
are between 1.1 and some high number.  It is believed
that, for each $d$, there is just one critical value $\lambda_d$
on $\Z^d$.  It is also to be expected that $\lambda_d$ is decreasing
in $d$, but only recently has it been shown that the largest critical
value on $\Z^d$ tends to~$0$ as $d \to \infty$.
This result was obtained by David Galvin and Jeff Kahn \cite{GK},
two combinatorialists, using graph theory, geometry, topology,
and lots of probabilistic combinatorics.  A consequence of their work
is that $\lambda=1$ is above the critical value(s) for sufficiently
large $d$; this can be stated in a purely combinatorial way: for
sufficiently high $d$ and large $n$, {\em most} independent sets in
$B_n^d$ are dominated by vertices of one parity.

In the next section we explain how we can use graphs to
understand models with hard constraints; then, in the section
following that, we will switch from $\Z^2$ to a much easier
setting, in which we can get our hands on nice Gibbs measures.

\section*{4. Hard constraints and graph homomorphisms} \addsec

\vskip-5mm \hspace{5mm}

We are interested in what are sometimes called ``nearest neighbor''
hard constraint models, where the constraints apply only to
adjacent sites.  Each site is to be assigned a ``spin'' from
some finite set, and only certain pairs of spins are permitted
on adjacent sites.  We can code up the constraints as a finite graph
$H$ whose nodes are the spins, and whose edges correspond to spins
allowed to appear at neighboring sites.  This {\em constraint graph}
$H$ may have some loops;
a loop at node $v \in H$ would mean that neighboring sites may
both be assigned spin $v$.  We adopt the statistical physics
tradition of reserving the letter ``$q$'' for the number of spins,
that is, the number of nodes in $H$.

The graph $G$ (e.g.\ $\Z^2$, above) of sites, usually infinite
but always countable and locally finite, is called (by us) the
{\em board}.  A legal assignment of spins to the sites of $G$
is nothing more or less than a graph homomorphism from $G$ to
$H$, i.e.\ a map from the sites of $G$ to the nodes of $H$ which
preserves edges.  We denote the set of homomorphisms from $G$
to $H$ by $\hom(G,H)$, and give it a graph structure by putting
$\vhi \sim \psi$ if $\vhi$ and $\psi$ differ at exactly one
site of $G$.

We will often confuse a graph with its set of nodes (or sites).
In particular, if $U$ is a subset of the nodes of $G$ then
$U$ together with the edges of $G$ contained in $U$ constitute
the ``subgraph of $G$ induced by $U$'', which we also denote by $U$.

In the hard-core model, the constraint graph $H$ consists of
two adjacent nodes, one of which is looped: a function from a
board $G$ to this $H$ is a homomorphism iff the set of
sites mapped to the unlooped node is an independent set.
Plate 4 shows some constraint graphs found in the literature.

When $H$ is complete and every node is looped as well, there
is no constraint and nothing interesting happens.

When $H$ is the complete graph $K_q$ (without loops), homomorphisms
to $H$ are just ordinary, ``proper'' $q$-colorings of the board.
(A proper $q$-coloring of a graph $G$ is a mapping from the nodes
of $G$ to a $q$-element set in which adjacent nodes are never
mapped to the same element.)
This corresponds to something called the ``anti-ferromagnetic
Potts model at zero temperature''.  In the anti-ferromagnetic
Potts model at positive temperature, adjacent sites are
merely discouraged (by an energy penalty), not forbidden,
from having the same spin; thus this is not a hard constraint
model in our terminology. The $q=2$ case of the Potts model is
the famous Ising model.

For a general constraint graph $H$, we need to elevate the notion
of activity to vector status.  To each node $i$ of $H$ we assign
a positive real activity $\lambda_i$, so that $H$ now gets an
activity vector $\lambda := (\lambda_1, \dots, \lambda_q)$.
When the board $G$ is finite, each homomorphism $\vhi \in
\hom(G,H)$ is assigned probability proportional to
$$
\prod_{v \in G} \lambda_{\vhi(v)}~.
$$
We can think of $\lambda_i$ as the degree to which we try to
use spin $i$, when it is available.  For example, if we know
the spins of the neighbors of site $v$ and consequently, say,
spins $i$, $j$ and $k$ are allowed for $v$, then the
$\lambda$-measure forces
$\Pr(\vhi(v)=i) = \lambda_i/(\lambda_i + \lambda_j + \lambda_k)$.

When $G$ is infinite, things get a little more complicated.
A finite subset (and its induced subgraph) $U \subset G$ will be
called a ``patch'' and its boundary $\partial U$ is the set
of sites not in $U$ but adjacent to some site of $U$.  We
define $U^+ := U \cup \partial U$.  If $\vhi$ is a function on $G$,
then $\vhi \restr U$ denotes its restriction to the subset $U$.

We say that $\mu$ is a Gibbs measure for $\lambda$ if:
for any patch $U \subset G$, and almost every $\psi\in \hom(G,H)$,
$$
\Pr_{\mu}\left(\vhi\restr U =\psi\restr U \given \vhi\restr (G-U)
=\psi\restr (G-U)\right) =
\Pr_{U^+}\left(\vhi \restr U =\psi\restr U \given \vhi\restr
\partial U =
\psi\restr \partial U \right)
$$
where ``$\Pr_{U^+}$'' refers to the finite $\lambda$-measure
on $U^+$.

This definition looks messy but it just means that
the probability distribution of a random $\vhi$ inside
a patch $U$ depends only on its value on the boundary of $U$, and is
the same as if $U$ and its boundary comprised all of the board.
We will see later that when $H$ has a certain nice property, as
it does in the case of the hard-core model, it suffices to check
the Gibbs condition only on patches consisting of a single
site---we call this the one-site condition.

It is a special case of a theorem of Dobrushin~\cite{D} that
there is always at least one Gibbs measure for any $\lambda$ on
$\hom(G,H)$; we are concerned with questions about when there is a
unique Gibbs measure, and when there is a phase transition (i.e. more
than one Gibbs measure).

Let us again look briefly at possible implications for phase transition
in the setting of finite boards.  Given a finite board $G$, a constraint
graph $H$ and activities $\lambda$, we define
the {\em point process} ${\cal P}(G,H,\lambda)$ as follows: starting from
any element of $\hom(G,H)$, choose a site $u$ of $G$ uniformly at random,
and give it a fresh spin according to the Gibbs condition, so that each
`legal' spin $j$ is chosen with probability proportional to $\lambda_j$.
The point process is a Markov chain on $\hom(G,H)$, and it is easy to check
that the $\lambda$-measure is a stationary distribution (which will be unique
provided $\hom(G,H)$ is connected, a point we will return to later).

Running the point process for sufficiently long will thus generate a random
homomorphism according to the $\lambda$-measure.  However, suppose that
the finite board $G$ is a large piece of an infinite board $G'$ exhibiting
a phase transition for our $\lambda$.  Then, if we start with a homomorphism
arising from one Gibbs measure on $\hom(G',H)$ (restricted to $G$), it is
reasonable to expect that the point process will take a long time to reach
a configuration resembling that from any other Gibbs measure on $\hom(G',H)$.
Thus it is generally believed that, in some necessarily loose sense,
phase transition on an infinite graph
corresponds to slow convergence for the point process on finite subgraphs.

\section*{5. Cayley trees and branching random walks} \addsec

\vskip-5mm \hspace{5mm}

Gibbs measures can be elusive and indeed it is generally a
difficult task to prove that phase transitions occur on a
typical board of interest, like $\Z^d$.  In order to get
results and intuition physicists sometimes turn to a more
tractable board, called by them the Bethe lattice (after
Hans Bethe) and by combinatorialists, usually, the Cayley tree.

We denote by $\T^r$ the $r$-branching Cayley tree, equivalently
the unique connected (infinite) graph which is cycle-free and
in which every site has degree $r\!+\!1$.  $\T^r$ is a vastly different
animal from $\Z^d$.  It is barely connected, falling apart
with the removal of any site; its patches have huge
boundaries, comparable in size with the patch itself;
its automorphism group is enormous.
It's surprising that we can learn anything at all about
$\hom(\Z^d,H)$ from $\hom(\T^r,H)$, and indeed we must
be careful about drawing even tentative conclusions in
either direction.  Basic physical parameters like entropy
become dodgy on non-amenable (big-boundary) boards like
$\T^r$ and a number of familiar statistical physics techniques
become useless.  More than making up for these losses, though,
are the combinatorial techniques we can use to study $\hom(\T^r,H)$.
There are even situations (e.g.\ in the study of information
dissemination) where $\T^r$ is the natural setting.

We are particularly interested in Gibbs measures on
$\hom(\T^r,H)$ which have the additional properties of being
{\em simple} and {\em invariant}.

For any site $u$ in
a tree $T$, let $d(u)$ be the number of edges incident with $u$ and
let $C_1(u), C_2(u), \dots, C_{d(u)}(u)$ be the connected
components of $T \setminus \{u\}$.

\begin{definition}\label{def:simple}
{\it A Gibbs measure $\mu$ on $\hom(T,H)$ is {\em simple} if, for
any site $u \in T$ and any node $i \in H$, the $\mu$-distributions
of
$$
\vhi\restr C_1(u), \dots, \vhi\restr C_{d(u)}(u)
$$
are mutually independent given $\vhi(u)=i$.}
\end{definition}

This condition, which is trivially satisfied by the
$\lambda$-measure for finite $T$, would follow
from the Gibbs condition itself if fewer than two of the
$C_i(u)$'s were infinite.

\begin{definition}\label{def:invariant}
{\it Let ${\mathcal A}(G)$ be the automorphism group of the board
$G$, and for any subset $S \subset \hom(G,H)$ and $\kappa \in
{\mathcal A}(G)$ let $S \circ \kappa := \{\vhi \circ \kappa:~\vhi
\in S\}$.  We say that a measure $\mu$ on $\hom(G,H)$ is {\em
invariant} if, for any $\mu$-measurable $S \subset \hom(G,H)$ and
any $\kappa \in {\mathcal A}(G)$, we have $\mu(S \circ \kappa) =
\mu(S)$.}
\end{definition}

Again, this condition is trivially satisfied for finite $G$;
but for an infinite board with as many automorphisms as $\T^r$, it
is quite strong.  Later we consider relaxing it slightly.
For now, we might well ask, how can we get our hands on any
Gibbs measure for $\hom(\T^r,H)$, let alone a simple, invariant one?

The absence of cycles in $\T^r$ makes it plausible that we
can get ourselves a Gibbs measure by building
random configurations in $\hom(\T^r,H)$ one
site at a time.  We could choose a root $x \in \T^r$, assign
it a random spin $i \in H$, then assign the neighbors of $i$
randomly to the $r\!+\!1$ children of $x$; thereafter, each
time a site $u$ gets spin $j$ we give its $r$ children random
spins from among the neighbors of $j$.

The process we have described can be thought of as a branching
random walk on $H$.  Imagine amoebas staggering from node
to adjacent node of $H$; each time an amoeba steps it divides into $r$
baby amoebas which then move independently at the next time step.
Of course, the (usually tiny) constraint graph $H$ is shortly
piled high with exponentially many amoebas, but being transparent
they happily ignore one another and go on stepping and dividing.

To get started we have to throw the first amoeba onto $H$
where we imagine that its impact will cause it to divide
$r\!+\!1$ ways instead of the usual $r$.

Note that the $r=1$ case is just ordinary random walk, started
somewhere on the doubly-infinite path $\T^1$ and run both
forward and backward.

To determine what probabilities are used in stepping from one
node of $H$ to an adjacent node, we assign a positive real
{\em weight} $w_i$ to each node.  For convenience we denote
by $z_i$ the sum of the weights of the neighbors of $i$ (including
$i$ itself, if there is a loop at $i$). An amoeba-child born on
node $i$ then steps to node $j$ with probability $w_j/z_i$.
If there is a loop at $i$, the amoeba stays at $i$ with
the appropriate probability, $w_i/z_i$.

Assuming $H$ is connected and not bipartite\footnote{A graph is
bipartite if its nodes can be partitioned into two sets neither
of which contains an edge.  Thus, for example, the existence
of a looped node already prevents $H$ from being bipartite.},
the random walk (branching or not) will have a stationary distribution
$\pi$; it is easily verified that $\pi_i$ is proportional to
$w_i z_i$ for each $i$, and somewhat less easily verified that
the mapping $w \hookrightarrow \pi$ is one-to-one provided
$\sum w_i$ has been normalized to 1.
We use the stationary distribution to pick the starting point for
the first amoeba, i.e.\ to assign a spin to the root of $\T^r$.

Finally, the payoff: not only does this node-weighted branching
random walk give us a simple invariant Gibbs measure; it's the
{\em only} way to get one.  The following theorem appears in
\cite{BW1} but it is not fundamentally different from characterizations
which can be found in Georgii \cite{G} and elsewhere.

\begin{theorem}\label{brw}
Let $H$ be a fixed connected
constraint graph with node-weights $w$ and let $r$ be
a positive integer.  Then the measure $\mu$ induced on $\hom(\T^r,H)$
by the $r$-branching $w$-random walk on $H$ is a simple,
invariant Gibbs measure, for some activity $\lambda$ on $H$.
Conversely, if $H$, $r$ and $\lambda$ are given, then every simple,
invariant Gibbs measure on $\hom(\T^r,H)$ is given by the $r$-branching
random walk on $H$ with nodes weighted by some $w$.
\end{theorem}

The proof is actually quite straightforward, and worth including here.
Invariance of $\mu$ with respect to root-preserving automorphisms of
$\T^r$ is trivial, since the random walk treats all children
equally; the only issue is whether the selection of root makes
a difference.  For this
we need only check that for two neighboring sites $u$ and $v$ of $\T^r$,
$\mu$ is the same whether $u$ is chosen as root or $v$ is.
But, either way we may choose $\vhi(u)$ and $\vhi(v)$ as
the first two spins and the rest of the procedure is the
same; so it suffices to check that for any (adjacent) nodes $i$ and $j$
of $H$, the probability that $\vhi(u)=i$ and $\vhi(v)=j$ is
the same with either root choice.  But these two probabilities
are
$$
\pi_i p_{ij} = z_iw_i \frac{w_j}{z_i} = w_iw_j = z_jw_j \frac{w_i}{z_j}
= \pi_j p_{ji}
$$
as desired.

To show that $\mu$ is simple is, indeed, simple: if we condition on
$\vhi(u)=i$ then, using invariance to put the root at $u$, the
independence of $\vhi$ on the $r\!+\!1$ components of $\T^r
\setminus \{u\}$ is evident from the definition of the branching
random walk.

The activity vector $\lambda$ for which $\mu$ is a Gibbs measure
turns out to be given by
$$
\lambda_i = \frac{w_i}{z_i^r}~.
$$

Let $U$ be any finite set of sites in $\T^r$, with
exterior boundary $\partial U$.  On account of invariance of labeling, we
may assume that the root $x$ does not lie in $U^+ = U \cup \partial U$.

Let $g \in \hom(U^+,H)$; we want to show that the probability that a
branching random walk $\vhi$ matches $g$ on $U$, given that it matches
on $\partial U$, is the same as the corresponding conditional probability
for the $\lambda$-measure.

Let $T$ be the subtree of $\T^r$ induced by $U^+$ and
the root $x$; for any $f \in \hom(T,H)$,
$$
\Pr(\vhi\restr T =f) = \pi_{f(x)} \cdot \prod_{u \rightarrow v}p_{f(u),f(v)}
$$
$$
= z_{f(x)}w_{f(x)} \prod_{u \rightarrow v}\frac{w_{f(v)}}{z_{f(u)}}
$$
where $u \rightarrow v$ means that $v$ is a child of $u$ in the tree.
The factors $z_{f(u)}$ corresponding to sites $u$ in $U$
each occur as denominator $r$ times in the above expression, since
each site in $U$ has all of its $r$ successors in $T$; and of course
each $w_{f(u)}$ occurs once as a numerator as well.  It follows that
if we compare $\Pr(\vhi\restr T =f)$ with $\Pr(\vhi\restr T =f')$,
where $f'$ differs
from $f$ only on $U$, then the value of the first is proportional to
$$
\prod_{u \in U}\frac{w_{f(u)}}{z_{f(u)}^r} = \prod_{u \in U}\lambda_{f(u)}
$$
which means that $\mu$ coincides with the finite measure, as desired.

Now let us assume that $\mu$ is a simple, invariant Gibbs measure on
$\hom(\T^r,H)$ with activity vector $\lambda$, with the intent of showing
that $\mu$ arises from a node-weighted branching random walk on $H$.

We start by constructing a $\mu$-random $\vhi$, site by site.
Choose a root $x$ of $\T^r$ and pick $\vhi(x)$
from the {\em a priori} distribution $\sigma$ of spins of $x$
(and therefore, by invariance, of any other site).  We next choose a
spin for the child $y$ of $x$ according to the conditional
distribution matrix $P = \{p_{ij}\}$ given by
$$
p_{ij} := \Pr\left(\vhi(y)=j \given \vhi(x)=i\right)~;
$$
again, by invariance of $\mu$, $P$ is the same for any pair of neighboring
sites.  It follows that $\sigma = \sigma \cdot P$, and moreover that
$P$ is the transition matrix of a reversible Markov chain, since the roles
of $x$ and $y$ can be interchanged.

Next we proceed to the rest of the children of $x$, then to the
grandchildren, etc., choosing each spin
conditionally according to all sites so far decided.

We claim, however, that the distribution of possible spins
of the non-root $v$ depends only on the spin of its parent $u$; this
is so because $\mu$ is simple and all sites so far ``spun'' are in
components of $\T^r \setminus \{u\}$ other than
the component containing $v$.  Thus the value of $\vhi(v)$ is given
by $P$ for {\em every} site $v \not= x$, and it follows that
$\mu$ arises from an $r$-branching Markov chain with state-space $H$,
starting at distribution $\sigma$.

Evidently for any (not necessarily distinct) nodes $i,~j$ of $H$, there
will be pairs $(u,v)$ of adjacent sites with $\vhi(u)=i$ and $\vhi(v)=j$
if and only if $i \sim j$ in $H$.  Hence $P$ allows transitions only
along edges of $H$, and there is a unique distribution $\pi$ satisfying
$\pi \cdot P = \pi$; thus $\sigma=\pi$.

It remains only to show that $P$ is a node-weighted random walk, and
it turns out that a special case of the Gibbs condition for one-site
patches suffices.
Let $j$ and $j'$ be nodes of $H$ which have a common neighbor $i$,
and suppose that all of the neighbors of the root $x$ have spin $i$.
Such a configuration will occur with positive probability and according
to the Gibbs condition for $U=\{x\}$,
$$
\frac{\Pr(\vhi(x)=j')}{\Pr(\vhi(x)=j)}=\frac{\lambda_{j'}}{\lambda_j}
$$
but
$$
\Pr(\vhi(x)=j)=\frac{\pi_j p_{ji}^{r+1}}{\sum_{k \sim i}\pi_k p_{ki}^{r+1}}
$$
and similarly for $j'$, so
$$
\frac{\Pr(\vhi(x)=j')}{\Pr(\vhi(x)=j)}
=\frac{\pi_{j'} p_{j'i}^{r+1}}{\pi_j p_{ji}^{r+1}}~.
$$
Thus the ratio
$$
\frac{p_{j'i}}{p_{ji}}= \left( \frac{\lambda_{j'} \pi_j}
{\lambda_j \pi_{j'}} \right)^{\frac{1}{r+1}}
$$
is independent of $i$.

Since $P$ is reversible we have $p_{ij} = \pi_j p_{ji} / \pi_i$, hence
$$
\frac{p_{ij'}}{p_{ij}} = \frac{\pi_{j'}p_{j'i}}{\pi_j p_{ji}}
$$
is also independent of $i$, and it follows that $P$ is a node-weighted
random walk on $H$.  This concludes the proof of Theorem~\ref{brw}.

\medskip

In view of Theorem~\ref{brw}, if we can understand the behavior of the
map $w \hookrightarrow \lambda$, we will know, given
$\lambda$, whether there is a nice Gibbs measure and if so
whether there is more than one.  The first issue is settled nicely
in the following theorem, a proof of which can be found in \cite{BW1}
and requires some topology.  A similar result was proved by Zachary~\cite{Z}.

\begin{theorem}
For every $r \ge 2$, every constraint graph $H$ and every set $\lambda$
of activities for $H$, there is a node-weighted branching random walk
on $\T^r$ which induces a simple, invariant Gibbs measure on $\hom(\T^r,H)$.
\end{theorem}

It's nice to know that we haven't required so much of our measures that
they can fail to exist.

A statistical physics dictum (true in great, but not unlimited, generality)
says that there's never a phase transition in dimension 1; that holds here:

\begin{theorem}\label{path}
For any connected constraint graph $H$ and any activity vector $\lambda$,
there is a unique simple invariant Gibbs measure on $\hom(\T^1,H)$.
\end{theorem}

Furthermore, in any dimension, there's always {\em some} region where
the map $w \hookrightarrow \lambda$ is one-to-one:

\begin{theorem}
For any $r$ and $H$ there is an activity vector $\lambda$ for
which there is only one simple invariant Gibbs measure on $\hom(\T^r,H)$.
\end{theorem}

\section*{6. Fertile and sterile graphs} \addsec

\vskip-5mm \hspace{5mm}

The fascination begins when we hit an $H$ and a $\lambda$ which
boast multiple simple, invariant Gibbs measures.  Let us examine
a particular case, involving a constraint graph we call the ``hinge''.

The hinge has three nodes, which we associate with the colors green, yellow
and red; all three nodes are looped and edges connect green with
yellow, and yellow with red.  Thus the only missing edge is
green-red, and a $\vhi \in \T^r$ may be thought of as a
green-yellow-red coloring of the tree in which no green site is
adjacent to a red one.

The hinge constraint in fact corresponds to a discrete version of the
Widom-Rowlinson model, in which two gases (whose particles
are represented by red and green) compete for space
and are not permitted to occupy adjacent sites; see e.g.\
\cite{BHW, WW, WR}.  When $\lambda_{\rm red}$ and $\lambda_{\rm green}$
are equal and large relative to $\lambda_{\rm yellow}$, the Widom-Rowlinson
model tends to undergo a phase transition as one gas spontaneously
dominates the other.  Plate 5 shows a red-dominated sample from the
Widom-Rowlinson model on $\Z^2$, with the unoccupied sites
left uncolored instead of being colored yellow.

We can see the phase transition operate on the Cayley tree $\T^2$.
If the green, yellow and red nodes are weighted 4, 2 and 1
respectively, $\lambda$ (normalized to integers) turns out to be
$(49, 18, 49)$---equal activity for green and red.  How can a
random walk which is biased so strongly toward green
end up coloring a tree according to a Gibbs measure with
symmetric specification?  As a clue, let us examine a site $u$
of $\T^2$ which happens to be surrounded by yellow neighbors.
To be colored green requires that a certain amoeba stepped from
yellow to green, then both of its children returned to yellow.
Thus the conditional probability that $u$ is green is proportional to
$$
\frac{4}{4+2+1}\cdot\left(\frac{2}{4+2}\right)^2
$$
as opposed to
$$
\frac{1}{4+2+1}\cdot\left(\frac{2}{2+1}\right)^2
$$
for red, but these values are equal.

Clearly the reversed weights 1, 2 and 4 would yield the same
activity vector, and in fact a third, symmetric weighting,
approximately 6, 7 and 6, does as well.  Plate 6 shows pieces
of $\T^2$ colored according to these three weightings.  Of
course the colorings have different proportions and are
easily identifiable; checking the stationary distributions
for the three random walks, we see that {\em a priori} a site
is colored green with probability about 59\% in the first
weighting, 30\% with the symmetric weighting and only 7\%
in the reversed weighting.  Yet, from a conditional point
of view, the three colorings are identical.

It turns out that the hinge is one of seven minimal graphs
each of which can produce a phase transition on $\T^r$ for
any $r \ge 2$.  The graphs are pictured in Plate 7.  We
say that a graph $H$ is {\em fertile} if $\hom(\T^r,H)$
has more than one simple, invariant Gibbs measure for
some $r$ and $\lambda$; otherwise it is {\em sterile}.
The fertile graphs are exactly those which contain one or
more of the seven baby graphs in Plate 7 as an induced
subgraph.  It turns out that the value of $r$ does not come into
play: if the constraint graph is rich enough to produce a phase
transition on any $\T^r$, then it does so for all $r\ge 2$.
One way to state the result is as follows:

\begin{theorem}\label{main}{\rm \cite{BW1}}
Fix $r>1$ and let $H$ be any constraint graph.  Suppose that
$H$ satisfies the following two conditions:
\medskip

{\rm (a)} Every looped node of $H$ is adjacent to all other nodes of $H$;
\medskip

{\rm (b)} With its loops deleted, $H$ is a complete multipartite graph.
\medskip

Then for every activity vector on $H$, there is a unique invariant
Gibbs measure on the space $\hom(\T^r,H)$.

If $H$ fails either condition {\rm (a)} or condition {\rm (b)} then
there is a set of activities $\lambda$ on $H$ for which $\hom(\T^r,H)$
has at least two simple, invariant Gibbs measures, and therefore $\lambda$
can be obtained by more than one branching random walk.
\end{theorem}

The proof of Theorem~\ref{main} is far too complex to reproduce
here, but reasonably straightforward in structure.  First, a
distinct pair of weightings yielding the same activity vector
must be produced for each of the seven baby fertile graphs, and
for each $r \ge 2$.  Second, it must be demonstrated that if
$H$ contains one of the seven as an induced subgraph, then
there are weightings (whose restrictions are close to those
previously found) which induce phase transitions on $\hom(\T^r,H)$.
Third, a monotonicity argument is employed to show that if $H$
satisfies conditions (a) and (b) of the theorem, then the map
from $w$ to $\lambda$ is injective.  Finally, an easy graph-theoretical
argument shows that $H$ satisfies (a) and (b) precisely if it
does not contain any of the seven baby fertile graphs as an
induced subgraph.

\section*{7. An application to statistical physics} \addsec

\vskip-5mm \hspace{5mm}

Theorem~\ref{main} has many shortcomings, applying as it does only
to hard constraint models on the Bethe lattice, and we must also
not forget that it considers only the very nicest Gibbs measures.
The constraint graph of the hard-core model is sterile, yet
it can have multiple Gibbs measures on $\T^r$ (or, as we saw,
on $\Z^2$) if we relax the invariance condition.

For $H$ = the hinge, however, and for any $r$ and $\lambda$, there
are multiple Gibbs measures on $\hom(\T^r,H)$ if and only if there
are multiple simple, invariant Gibbs measures.  Like the Ising
model, the Widom-Rowlinson model exhibits spontaneous symmetry-breaking;
indeed its relationship to the Ising model parallels the relation
between nodes and edges of a graph.

One of the nice properties
known for the Ising model is that it can exhibit at most one
critical point; but the proof of this fact does not work for the
Widom-Rowlinson model.  Indeed, in \cite{BHW} the methods above
are used to construct a board $G$ for which $\hom(G,H)$ has
three (or more) calculable critical points, with $H$ the hinge.
Set $\lambda = \lambda_{\rm green} = \lambda_{\rm red}$, fixing
$\lambda_{\rm yellow}=1$, so that the single parameter $\lambda$
controls the Widom-Rowlinson model.  Then:

\begin{theorem}\label{thm:WR}
There exist $0<\lambda_1<\lambda_2<\lambda_3$ and an infinite graph $G$,
such that the Widom--Rowlinson model on $G$ with activity $\lambda$ has a
unique Gibbs measure for
$\lambda\in (0,\lambda_1] \cup [\lambda_2, \lambda_3]$,
and multiple Gibbs measures for
$\lambda\in (\lambda_1,\lambda_2)\cup (\lambda_3,\infty)$.
\end{theorem}

The board $G$ constructed in \cite{BHW} is a tree, but not quite
a regular one; it is made by dangling seven new pendant sites from
each site of $\T^{40}$.  Readers are referred to that paper for
the calculations, but the intuition is something like this.

For low $\lambda$ the random coloring
of $G$ is mostly yellow, but as $\lambda$ rises, either green or
red tends to take over the interior vertices as in $\T^{40}$.
Then comes the third interval, where the septuplets of leaves, wanting
to use both green and red, force more yellow on the interior
vertices, relieving the pressure and restoring green-red symmetry.
Finally the activity becomes so large that the random coloring
is willing to give up red-green variety among the septuplets in
order to avoid yellow interior vertices, and symmetry-breaking
appears once again.

It turns out that multiple critical points can be obtained for
the hard-core model in a similar way.

\section*{8. Dismantlable graphs} \addsec

\vskip-5mm \hspace{5mm}

In addition to the fertile and sterile graphs, a second graph
dichotomy appears repeatedly in our studies: dismantlable and
non-dismantlable graphs.  Coincidentally, the term ``dismantlable''
as applied to graphs was coined by Richard Nowakowski and the second
author \cite{NW} almost twenty years ago in another context
entirely: a pursuit game on graphs.

Two players, a cop ${\mathcal C}$ and a robber ${\mathcal R}$,
compete on a fixed, finite, undirected graph $H$.  We will assume
that $H$ is connected and has at least one edge, although the concepts
make sense even without these assumptions.
The cop begins by placing herself at a node of her choice; the
robber then does the same.  Then the players alternate
beginning with ${\mathcal C}$, each moving to an adjacent
node.  The cop wins if she can ``capture'' the robber,
that is, move onto the node occupied by the robber;
${\mathcal R}$ wins by avoiding capture indefinitely.
In doing so ${\mathcal R}$ is free to move (or even place
himself initially) onto the same node as the cop, although
that would be unwise if the node were looped since then
${\mathcal C}$ could capture him at her next move.

Evidently the robber can win on any loopless graph by
placing himself at the same node as the cop and then shadowing
her every move; among graphs in which every node is looped,
${\mathcal C}$ clearly wins on paths and loses on cycles of length
4 or more.  (In the game as defined in \cite{NW,Q}, there is in
effect a loop at every node of $H$.)

The graph on which the game is played is said to be {\em cop-win} if
${\mathcal C}$ has a winning strategy, {\em robber-win} otherwise.
The following structural characterization of cop-win graphs is proved
in \cite{NW} for the all-loops case, but in fact the proof (which is
not difficult, and left here as an exercise) works fine in our more
general context.

Let $N(i)$ be the neighborhood of node $i$ in $H$ and suppose
there are nodes $i$ and $j$ in $H$ such that $N(i) \subseteq N(j)$.
Then the map taking $i$ to $j$, and every other node of $H$ to itself,
is a homomorphism from $H$ to $H \setminus \{i\}$.
We call this a {\em fold} of the graph $H$.
A finite graph $H$ is {\em dismantlable} if there is a sequence of folds
reducing $H$ to a graph with one node (which will necessarily be looped).

Note that dismantlable graphs are easily recognized in polynomial time.
Plate 8 shows some dismantlable and non-dismantlable graphs.

The following theorem, from \cite{BW2}, collects a boatload of equivalent
conditions.

\begin{theorem}\label{equivalent}
The following are equivalent, for finite connected graphs $H$ with at least
one edge.
\begin{enumerate}
\item $H$ is dismantlable.
\item $H$ is cop-win.
\item For every finite board $G$, $\hom(G,H)$ is connected.
\item For every board $G$, and every pair $\vhi, \psi \in \hom(G,H)$
agreeing on all but finitely many sites, there is a path in $\hom(G,H)$
between $\vhi$ and $\psi$.
\item There is some positive integer $m$ such that, for every board $G$,
every pair of sets $U$ and $V$ in $G$ at distance at least $m$,
and every pair of
maps $\vhi, \psi \in \hom(G,H)$, there is a map $\theta \in \hom(G,H)$
such that $\theta$ agrees with $\vhi$ on $U$ and with $\psi$ on $V$.
\item For every positive integer $r$, and every pair of maps $\vhi, \psi \in
\hom(\T^r,H)$, there is a site $u$ in $\T^r$ with $\vhi(u) \not= \psi(u)$,
a patch $U$ containing $u$, and a map $\theta \in \hom(\T^r,H)$ which agrees
with $\psi$ on $\T^r \setminus U$ and with $\vhi$ on $u$.
\item For every board $G$ and activity vector $\lambda$, if $\mu$ is a
measure on $\hom(G,H)$ satisfying the one-site condition, then $\mu$ is a
Gibbs measure.
\item For every finite board $G$ and activity vector $\lambda$,
every stationary distribution for the point process ${\cal P}(G,H,\lambda)$
is a Gibbs measure.
\item For every board $G$ of bounded degree such that $\hom(G,H)$ is
non-empty, there is an activity vector $\lambda$ such that
there is a unique Gibbs measure on $\hom(G,H)$.
\item For every $r$, there is an
activity vector $\lambda$ such that there is a unique Gibbs measure on
$\hom(\T^r,H)$.
\end{enumerate}
\end{theorem}

We will prove here what we think, to a graph theorist, is the most
interesting of these equivalences---$(i)$ and $(iii)$.  Recall that
two maps in $\hom(G,H)$ are adjacent if they differ on one site
of $G$.

Let us first assume $H$ is dismantlable.  If it has only one
node the connectivity of $\hom(G,H)$ is trivial, since it has
at most one element.  Otherwise there are nodes $i \not= j$ in
$H$ with $N(i) \subseteq N(j)$ and we may assume by induction
that $\hom(G,H')$ is connected for $H' := H \setminus \{i\}$.

Define, for $\vhi$ in $\hom(G,H)$, the map $\vhi'$ in
$\hom(G,H')$ (and also in $\hom(G,H)$) by changing all $i$'s to
$j$'s in the image.  If $\alpha$ and $\beta$ are two maps in
$\hom(G,H)$ then there are paths from $\alpha$ to $\alpha'$,
$\alpha'$ to $\beta'$ and $\beta'$ to $\beta$; so $\hom(G,H)$
is connected as claimed.

For the converse, let $H$ be non-dismantlable, and suppose that
nonetheless $\hom(G,H)$ is connected for all finite boards $G$;
let $q = |H|$ be minimal with respect to these properties.

If there are nodes $i$ and $j$ of $H$ with $N(i) \subseteq N(j)$, then
$H := H \setminus \{i\}$ is also non-dismantlable. In this case, we claim
that the connectivity of $\hom(G,H)$ implies connectivity of
$\hom(G,H')$.  To see this, define, for $\vhi \in \hom(G,H)$,
the map $\vhi' \in \hom(G,H')$ by changing all $i$'s to $j$'s
in the image as before.  If $\alpha$ and $\beta$ are two maps in
$\hom(G,H')$, then we may connect them by a path
$\vhi_1, \dots, \vhi_t$ in $\hom(G,H)$; now we observe that
the not-necessarily distinct sequence of maps $\vhi'_1, \dots,
\vhi'_t$ connects $\alpha$ and $\beta$ in $\hom(G,H')$.
This contradicts the minimality of $H$, so we may assume from now on that
there is no pair of nodes $i \not= j$ in $H$ with $N(i) \subseteq N(j)$.

Now let $G$ be
the `weak' square of $H$, that is, the graph whose
nodes are ordered pairs $(i_1,i_2)$ of nodes of $H$ with
$(i_1,i_2) \sim (j_1,j_2)$ just when $i_1 \sim j_1$ and
$i_2 \sim j_2$.  There are two natural homomorphisms from
$G$ to $H$, the projections $\pi_1$ and $\pi_2$, where
$\pi_1(i_1,i_2)=i_1$ and $\pi_2(i_1,i_2)=i_2$; we claim
that $\pi_1$ is an isolated point of the graph $\hom(G,H)$, which
certainly implies that $\hom(G,H)$ is disconnected.

If not, there is a map $\pi'$ taking (say) $(i_1,i_2)$ to
$k \not= i_1$ and otherwise agreeing with $\pi_1$.  Let $j_2$
be a fixed neighbor of $i_2$ and $j_1$ any neighbor of $i_1$.
Then $(i_1,i_2) \sim (j_1,j_2)$, and hence $k \sim j_1$ in $H$.
We have shown that every neighbor of $i_1$ is also a neighbor of $k$,
contradicting the assumption that no such pair of nodes exists in $H$.
This completes the proof.

For the last part of the proof, there is also a simpler (and smaller)
construction that works provided $H$ has at least one loop: see~\cite{BW2}
or Cooper, Dyer and Frieze~\cite{CDF}.

We have seen that, for a dismantlable constraint graph $H$, and any
board $G$ of bounded degree, there is some $\lambda$ (which can be
taken to depend only on $H$ and the maximum degree of $G$) such that
there is a unique Gibbs measure on $\hom(G,H)$.  Dyer, Jerrum and
Vigoda~\cite{DJV} have proved a ``rapid mixing'' counterpart to this
result: given a dismantlable $H$, and a degree bound $\Delta$, there
is some $\lambda$ such that the point process ${\cal P}(G,H,\lambda)$
is rapidly mixing for all finite graphs $G$ with maximum degree at most
$\Delta$.  Of course, if $H$ is not dismantlable, then no such result
can be true as $\hom(G,H)$ need not be connected.

\section*{9. Random colorings of the cayley tree} \addsec

\vskip-5mm \hspace{5mm}

We have observed that ordinary ``proper'' $q$-colorings of a graph
$G$ are maps in $\hom(G,K_q)$; since $K_q$ is sterile, there
is never more than one Gibbs measure for $G = \T^r$.  However,
we see even in the case $q=2$ that multiple Gibbs measures
exist, because each of the two 2-colorings of $\T^r$ determines
{\em by itself} a trivial Gibbs measure, as does any convex
combination.  However, only the $\frac12$, $\frac12$ combination
is invariant under parity-changing automorphisms of $\T^r$.

All of the Gibbs measures in the $q=2$ case are, however, simple and
invariant under all the parity-preserving automorphisms of the
board.  Such Gibbs measures are neededed to realize the phase transition
for the hard-core model as well, so it is not surprising that it
is useful to relax our requirements slightly and to consider these
{\em semi-invariant} simple Gibbs measures.

Fortunately we don't have to throw away all our work on invariant
Gibbs measures in moving to semi-invariant ones.  Given a
constraint graph $H$ on nodes $1,2, \dots, q$ which is connected
and not bipartite, we form its bipartite ``double'', denoted $2H$,
as follows: the nodes of $2H$ are $\{1,2\dots,q\} \cup
\{-1, -2, \dots, -q\}$ with an edge between $i$ and $j$ just when
$i \sim -j$ or $-i \sim j$ in $H$.  Note that $2H$ is loopless;
a loop at node $i$ in $H$ becomes the edge $\{-i,i\}$ in $2H$.

A homomorphism $\psi$ from $\T^r$ to $2H$ induces a homomorphism
$|\psi|$ to $H$ via $|\psi|(v) = |\psi(v)|$.  In the reverse
direction, a map $\vhi$ in $\hom(\T^r,H)$ may be transformed to
a map $\bar\vhi$ in $\hom(\T^r,2H)$,
by putting $\bar\vhi(v)=\vhi(v)$ for even sites $v \in \T^r$
and $\bar\vhi(v)=-\vhi(v)$ for odd $v$.

Let $\lambda = (\lambda_1, \dots, \lambda_n)$ be an
activity vector for $H$ and suppose that $\mu$ is a simple invariant Gibbs
measure on $\hom(\T^r,H)$ corresponding to $\lambda$.  From $\mu$
we can obtain a simple invariant Gibbs measure $\bar\mu$ on
$\hom(\T^r,2H)$ by selecting $\vhi$ from $\mu$, and flipping a fair
coin to decide between $\bar\vhi$ (as defined above) and $-\bar\vhi$.
Obviously $\bar\mu$ yields the activity vector $\bar \lambda$ on $2H$
given by $\bar\lambda_i = \lambda_{|i|}$.  Furthermore, the weights on
$H$ which produce $\mu$ extend to $2H$ by $w_{-i}=w_i$.

Conversely, suppose $\nu$ is a simple invariant Gibbs measure on
$\hom(\T^r,2H)$ whose activity vector satisfies $\lambda_{-i}=\lambda_i$
for each $i$.  Then the measure $|\nu|$, obtained by choosing $\psi$
from $\nu$ and taking its absolute value, is certainly an invariant
Gibbs measure on $\hom(\T^r,H)$ for $\lambda \restr \{1, \dots, q\}$,
but is it simple?

In fact, if the weights on $2H$ which produce $\nu$ do not satisfy
$w_{-i}=cw_i$, then $|\nu|$ will fail to be simple.  To see this,
observe that if the weights are not proportional then there are nodes
$i \sim j$ of $H$ such that $p_{-i,-j} \not= p_{i,j}$ in the random walk
on $2H$.  Suppose that $|\psi|$ is conditioned on the color of the root
$w$ of $\T^r$ being fixed at $i$, and let $x$ and $y$ be distinct
neighbors of $w$. Set $\alpha = \Pr(\psi(w)=i \mid |\psi(w)| = i)$. Then
$$
\Pr(|\psi|(x)=j) = (1-\alpha) p_{-i,-j} + \alpha p_{i,j}
$$
but
$$
\Pr(|\psi|(x)=j \wedge |\psi|(y)=j) = (1-\alpha) p_{-i,-j}^2 +
\alpha p_{i,j}^2 > \Pr(|\psi|(x)=j)^2
$$
so the colors of $x$ and $y$ are not independent given $|\psi|(w)$.

However, we can recover simplicity at the expense of one bit worth of
symmetry.  Let $\nu^+$ be $\nu$ conditioned on $\psi(u)>0$, and
define $\nu^-$ similarly.  Then $|\nu^+|$ and $|\nu^-|$ are essentially
the same as $\nu^+$ and $\nu^-$, respectively, and all are simple; but
these measures are only semi-invariant.

On the other hand, suppose $\mu$ is a simple semi-invariant Gibbs
measure on $\hom(\T^r,H)$.  Let $\theta$ be a parity-reversing
automorphism of $\T^r$ and define $\mu' := \mu \circ \theta$,
so that $\frac12 \mu + \frac12 \mu'$ is fully invariant (but generally
no longer simple).  However, $\nu := \frac12 \bar\mu + \frac12
(-\overline{\mu'})$
is a simple {\em and} invariant Gibbs measure on $\hom(\T^r,2H)$, thus
given by a node-weighted random walk on $2H$.  We can recover $\mu$
as $\nu_+$, hence:

\begin{theorem}\label{thm:semi}
Every simple semi-invariant Gibbs measure on $\hom(\T^r,H)$ is obtainable
from a node-weighted branching random walk on $2H$, with its initial state
drawn from the stationary distribution on positive nodes of $2H$.
\end{theorem}

Suppose, instead of beginning with a measure, we start by weighting the
nodes of $2H$ and creating a Gibbs measure as in Theorem~\ref{thm:semi}.
Suppose the activities of the measure are $\{\lambda_i:~i=\pm 1, \dots,
\pm q\}$.  By identifying color $-i$ with $i$ for each $i>0$, we create
a measure on $H$-colorings, but this will not be a Gibbs measure unless
it happens that $(\lambda_{-1}, \dots, \lambda_{-q})$ is proportional
to $(\lambda_1, \dots, \lambda_q)$.

We could assure this easily enough by making the weights proportional
as well, e.g.\ by $w_{-i}=w_i$; then the resulting measure on $\hom(G,H)$
could have been obtained directly by applying these weights to $H$,
and is thus a fully invariant simple Gibbs measure.  To get new,
semi-invariant Gibbs measures on $\hom(G,H)$, we must somehow devise
weights for $2H$ such that $w_{-i} \not\propto w_i$ yet
$\lambda_{-i} \propto \lambda_i$.

Restated with slightly different notation, simple semi-invariant Gibbs
measures are in 1--1 correspondence with solutions to the ``fundamental
equations''
$$
\lambda_{i} = \frac{u_i}{\left( \sum_{j \sim i}v_j \right) ^r}
= \frac{v_i}{\left( \sum_{j \sim i}u_j \right) ^r}
$$
for $i = 1, \dots, q$.  Such a solution will be invariant if
$u_i = v_i$ for each $i$.

Plate 9 illustrates a semi-invariant, but not
invariant, simple Gibbs measure for uniform 3-colorings of $\T^3$.
Approximate weights of the nodes of $2H = 2K_3$ are given along with
part of a sample coloring drawn from this measure.  Additional
measures may be obtained by permuting the colors or by making all
the weights equal (invariant case).

Results for $q$-colorings of $\T^r$, with $q>2$ and $r>1$, are
as follows:

When $q < r\!+\!1$,  all choices of activity vector including
the uniform case yield multiple simple semi-invariant Gibbs measures.

When $q > r\!+\!1$, there is only one simple semi-invariant Gibbs
measure for the uniform activity vector, but multiple simple
semi-invariant Gibbs measures for some other choices of activity vector.

The critical case is at $q=r\!+\!1$, that is, when the number of
colors is equal to the degree of the Cayley tree.  Here it turns out that
there are multiple simple semi-invariant Gibbs measures for all activity
vectors {\em except} the uniform case, where there is just one.

When $q > r\!+\!1$ and the activities are equal, the
unique simple semi-invariant Gibbs measure is in fact the only Gibbs measure
of any kind.  This was conjectured in \cite{BW3} but proved only for
$q>cr$, with fixed $c>1$; Jonasson \cite{J} has recently, and
very nicely, finished the job.  Jonasson's result is in a sense a special
case of the conjecture that the Markov chain of $q$-colorings of a
finite graph of maximum degree less than $q$, which progresses by choosing
and recoloring sites randomly one at a time, mixes rapidly.  So far
the best result is Vigoda's \cite{V} which proves this if the maximum
degree is at most $6q/11$.

When $q \le r\!+\!1$ there are lots of other Gibbs measures, including
ones we call {\em frozen}.  These come about because it is
possible for a measure to satisfy the Gibbs condition in
a trivial and somewhat unsatisfactory way.  For example, suppose
we are $q$-coloring $\T^r$ (with root $w$) for some $q \le r\!+\!1$,
and let $\psi$ be any fixed coloring in which the children of every node
exhibit all colors other than the color of the parent.  Let $\mu$ be
the measure which assigns probability 1 to $\psi$.  Then for any finite
patch $U$, which we can assume to be a subtree including the root, the colors
on $\partial U$ force the colors on the leaves of $U$, and we can continue
inwards to show that the original coloring $\psi \restr U$ is the only one
consistent with the colors on $\partial U$.
Thus $\mu$ satisfies the Gibbs condition trivially, and is
also vacuously simple---but not invariant or semi-invariant.  We call
a Gibbs measure of this type ``frozen''.
A frozen state of $\hom(\T^2,K_3)$ is illustrated in Plate 10.
For more about frozen Gibbs measures the reader is referred to \cite{BW2}.

In a soft constraint model such as the Potts model, frozen Gibbs measures
can only occur at zero temperature.  Since most of the time statistical
physicists are interested only in phases which exist at some positive
temperature (and have positive entropy), frozen measures are generally
absent from the statistical physics literature.  However, they are
interesting combinatorially and motivate some definitions in the
next section.

\section*{10. From statistical physics back to graph theory} \addsec

\vskip-5mm \hspace{5mm}

We conclude these notes with a theorem and a conjecture in ``pure''
graph theory, stripped of probability and physics, but suggested by
the many ideas which have appeared in earlier sections.

Suppose $H$ is bipartite and we are given some $\vhi \in \hom(\T^1,H)$,
where the sites of $\T^1$ are labeled by the integers $\Z$.  Then
knowing $\vhi(n)$ even for a very large $n$ tells us something about
$\vhi(0)$, namely which ``part'' of $H$ it is in.  We call this
phenomenon {\em long range action}, and define it on Cayley
trees as follows:  If there is a $\vhi \in \hom(\T^r,H)$ and a
node $i \in H$ such that for any $n$, no $\psi \in \hom(\T^r,H)$
agreeing with $\vhi$ on the sites at distance $n$ from the root
can have spin $i$ at the root, we say $\hom(\T^r,H)$ has long
range action.

\begin{theorem}\label{warm}
If $H$ is $k$-colorable then $\hom(\T^{k-1},H)$ has long range action.
\end{theorem}

For example, the coloring described at the end of the previous section, which gives rise to a frozen Gibbs
measure, shows that $\hom(\T^2,K_3)$ (more generally, $\hom(\T^r,K_{r+1})$) has long range action.  We also see
from Theorem~\ref{equivalent} that \linebreak $\hom(\T^r,H)$ has long range action for {\em no} $r$ if and only if
$H$ is dismantlable; of course then $H$ has at least one loop and therefore has infinite chromatic number.

Note that the theorem connects a statement about homomorphisms
{\em from} $H$ to a statement about homomorphisms {\em to} $H$.
However, it is difficult to see how to turn a $k$-coloring of
$H$ into a suitable map in $\hom(\T^{k-1},H)$.  Suppose, for
instance, that $H$ is the 5-cycle $C_5$, with nodes represented
by the integers modulo 5.  We can get a
completely frozen map in $\hom(\T^2,C_5)$ by making sure we
use both $i\!+\!1$ and $i\!-\!1$ on the children of any site
of spin $i$.  But what has this map got to do with any 3-coloring
of $C_5$?

The proof of Theorem~\ref{warm}, found in \cite{BW4}, uses a
vector-valued generalization of coloring to construct the
required map in $\hom(\T^{k-1},H)$.

We now move from long range action to the familiar notion
of connectivity.  Theorem~\ref{equivalent}---in fact, the part
whose proof is given above---tells us that $\hom(G,H)$ is
connected for any finite $G$ just when $H$ is dismantlable.
Suppose we restrict ourselves to boards of bounded degree?
If, for example, $H$ is bipartite, $\hom(K_2,H)$ is already
disconnected.  If $H=K_q$ then $\hom(K_q,H)$ is extremely
disconnected, consisting of $n$!\ isolated maps.
By analogy with Theorem~\ref{warm}, we should perhaps be able to prove:

\begin{conjecture}\label{mobile}
If $H$ is $k$-colorable then $\hom(G,H)$ is disconnected for
some finite $G$ of maximum degree less than $k$.
\end{conjecture}

A proof for $k=3$ appears in \cite{BW4} and Lov\'asz \cite{L}
has shown that the conjecture holds for $k=4$ as well.  We
think that a proof of Conjecture~\ref{mobile} would have to
capture some basic truths about graphs and combinatorial
topology, and fervently hope that some reader of these notes
will take up the challenge.

\label{lastpage}

\end{document}